\documentstyle[12pt]{article}

  \newcommand{\const}{\rm const}

\begin{document}

 \begin{center}

 \ {\bf  Analog of modulus of convexity for Grand Lebesgue Spaces.} \\

  \vspace{6mm}

 {\bf M.R.Formica, \ E.Ostrovsky, L.Sirota. }

 \end{center}

\vspace{6mm}

 \ Universit\`{a} degli Studi di Napoli Parthenope, via Generale Parisi 13, Palazzo Pacanowsky, 80132,
Napoli, Italy. \\

e-mail: mara.formica@uniparthenope.it \\

\vspace{5mm}

 \ Israel,  Bar - Ilan University, department  of Mathematic and Statistics, 59200, \\

\vspace{4mm}

e-mails: \  eugostrovsky@list.ru \\

\vspace{4mm}

\begin{center}

  {\bf Abstract} \\

\end{center}

 \ We introduce and evaluate the degree of convexity of an unit ball, so - called, characteristic of convexity  (COC) for
 the Grand Lebesgue Spaces, (GLS), which is a slight analog of the classical notion of the modulus of convexity (MOC).

\vspace{5mm}

 \ {\it Key words and phrases.} Banach, Lebesgue - Riesz and Grand Lebesgue Spaces (GLS) and norms, triangle inequality, unit
 ball, embedding, modulus of convexity (MOC), weak characteristic of convexity  (WCOC).

\vspace{5mm}

\section{Notations. Definitions. Statement of problem.}

\vspace{4mm}

 \ Let $ \ (X,||\cdot||) \ $ be Banach space, $ \ S \ $ be its unit sphere: $ \ S = \{ \ x, \ x \in X, \ ||x|| = 1 \ \}  \ $
and $ \ B \ $ be its unit ball with the center in origin: $ \ B = \{ \ x, \ x \in X, ||x|| \le 1.  \ \}  $  Let also $ \ \epsilon \ $
be arbitrary number from the segment $ \  [0,2]: \ 0 \le \epsilon \le 2. \ $ Recall the following
very important in the geometrical theory of Banach spaces notion Modulus Of Convexity  (MOC)  $ \ \delta_X(\epsilon) \ $  for the space
$ \ (X,||\cdot||): \ $

\vspace{3mm}

\ {\bf Definition 1.1.} The Modulus Of Convexity  (MOC) for the space \\
 $ \ X = (X,||\cdot||) = (X,||\cdot||X), \ $ which is denoted by $ \ \delta_X(\epsilon) \ $ is defined as follows

\begin{equation} \label{MOC Ball}
\delta_X(\epsilon)  \stackrel{def}{=} \inf \left\{ \ 1 - \frac{||x + y||}{2}: \ x,y \in B; \ ||x - y|| \ge \epsilon  \ \right\};
\end{equation}
the {\it Ball definition}; or equally

\begin{equation} \label{MOC sphere}
\delta_X(\epsilon)  \stackrel{def}{=} \inf \left\{ \ 1 - \frac{||x + y||}{2}; \ x,y \in S; \ ||x - y|| \ge \epsilon  \ \right\};
\end{equation}
spherical definition.\par

\vspace{4mm}

\ This important  for Functional Analysis notion was introduced by  O.Hanner (1956), see  \cite{Hanner}, and was investigated in many works,  see e.g.
\cite{Clarkson}, \cite{FigeL1}, \cite{FigeL2}, \cite{Gao},  \cite{Hanner}, \cite{Pick Kufner John Fucik}, \cite{Reif} and so one. \par

 \ The other application, indeed, in the theory of random fields, may be found in \cite{Ostrovsky1}, chapter 3. \par

\vspace{3mm}

 \ For example, let $ \ X \ $ be the classical Lebesgue - Riesz space $ \ L_p, \ $ builded on certain atomless measure space;
 the correspondent Module Of Convexity will be denoted by $ \ \delta_p(\epsilon);  \ p \in (1,\infty). \ $  If $ \  p \in (1,2], \ $ then
the function $ \ \delta_p(\epsilon) \ $ is an unique positive solution of an equation

$$
\left( \  1 - \delta_p(\epsilon) + 0.5 \epsilon  \ \right)^p  + \left( \  1 - \delta_p(\epsilon) - 0.5 \epsilon  \ \right)^p = 2,
$$
so that when $ \ \epsilon \in [0,2]  \ $

\begin{equation} \label{small p}
\delta_p(\epsilon) \ge \frac{p-1}{8} \ \epsilon^2.
\end{equation}

 \vspace{3mm}

 \ If now $ \ p \in (2,\infty), \ $  then

\begin{equation} \label{great p beg}
\delta_p(\epsilon) = 1 - \left( \ 1 - (0.5 \epsilon)^p  \ \right)^{1/p},
\end{equation}
and when again  $ \ \epsilon \in [0,2] \ $

\begin{equation} \label{great p appr}
\delta_p(\epsilon) \ge \frac{\epsilon^p}{ p \ 2^p};
\end{equation}
see e.g. \cite{FigeL1}, \cite{FigeL2}, \cite{Gao}.\par

\vspace{4mm}

 \ {\bf  We intent in this short report to introduce some weak analog of modulus of convexity for
 Grand Lebesgue Spaces (GLS) and derive some its properties. } \par

\vspace{3mm}

 \ It follows from the definition 1.1. that for $ \  x,y \in B  \ $

\begin{equation} \label{key ineq}
||x + y|| \le 2 - 2 \delta_X (||x - y||),
\end{equation}

a refined triangle inequality. \par

\vspace{3mm}

 \ Leu us give some generalization of this notion MOC. \par

\vspace{4mm}

 \ {\bf Definition 1.2.} Let again $ \ (X,||\cdot||X = ||\cdot|| ) \ $ be the Banach space. Suppose that there exists
 an {\it another} Banach space $ \ (Y,|||\cdot|||Y = |||\cdot||| ) \ $ such that $ \ X \ $ is embedded in $ \ Y: \  X \subset Y, \ $
  and  a  non - negative numerical valued  function  (functional!)  $ \  \ \Delta[X,Y] = \Delta[X](u), \ u \in X, \ $  which is named
 as a {\it weak characteristic of convexity,  (WCOC),}  such that

$$
\Delta[X](u) = \Delta[X,Y](u) = 0 \ \Leftrightarrow u = 0,
$$
 and for which

\begin{equation} \label{weak key notion}
\forall x,y \in B \ \Rightarrow ||x + y|| \le 2 - 2 \Delta[X,Y](|||x - y|||).
\end{equation}

\vspace{4mm}

 \ It is this inequality  (\ref{weak key notion})    for the Grand Lebesgue Spaces $ \ X, \ $
 that was applied in particular in the theory if random fields,  see for example \cite{Ostrovsky1}, chapter 3, sections 3.3 - 3.6. \par

\vspace{4mm}

 \ To be more precisely, we want to prove the existence  $ \ \Delta[G\psi](\cdot) \ $ for the so - called Grand Lebesgue Spaces $  \ G\psi \ $
 and  derive  some its estimations. \par

\vspace{4mm}

\begin{center}

\ {\sc Brief note about  Grand Lebesgue Spaces (GLS). } \\

\end{center}

\vspace{5mm}

 \ We recall here for reader convenience some known definitions and  facts about  the theory of Grand Lebesgue Spaces (GLS) using in this article.
Let $ \ (Z, M, \mu) \ $ be measurable space with non - trivial atomless measure $ \ \mu. \ $ The ordinary Lebesgue - Riesz norm $ \ ||f ||_p \ $
for the numerical valued measurable function $ \ f: Z \to R \ $ is defined as ordinary

$$
||f||_p := \left[ \  \int_Z |f(z)|^p \ \mu(dz)  \ \right]^{1/p}, \ 1 \le p < \infty;
$$
and as ordinary $ \ L_p  = L_p(Z,\mu) = \{ \ f: Z \to R, \ ||f||_p < \infty  \ \}. \ $ \par

\vspace{3mm}

   \ Further, let  the numbers $ \  (a,b) \ $ be constants  such that $ \ 1 \le a < b \le \infty; \ $   and let
   $ \ \psi = \psi(p) = \psi[a,b](p), \ p \in (a,b), \ $ be numerical valued strictly positive function  not necessary to be finite in every point:

\begin{equation} \label{Positive psi}
 \inf_{p \in (a,b)} \psi[a,b](p) > 0.
\end{equation}

\vspace{3mm}

 \ In the case when $ \ b < \infty \ $ one can assume sometimes $ \ p \in (a,b]. \ $  \par

  \    For instance

  $$
    \psi_{(m)}(p) := p^{1/m}, \ m = \const > 0, \ p \in [1,\infty)
  $$
or

$$
   \psi^{(b; \beta)}(p) := (p-a)^{-\beta_1} \cdot (b-p)^{-\beta_2}, \ p \in (a,b), \ \beta_1, \ \beta_2 = \const \ge 0.
$$
 \ The set of all such a functions will be denoted by $ \ \Psi = \Psi[a,b] = \{  \psi(\cdot)  \}. \ $ \par

 \ By definition, the (Banach) Grand Lebesgue Space (GLS)    $  \ G \psi  = G\psi [a,b],  $
    consists on all the real (or complex) numerical valued measurable functions
   $   \  f: Z \to R \ $  defined on the whole our  space $ \ Z \ $ and having a finite norm

 \begin{equation} \label{norm psi}
    || \ f \ ||G\psi = ||f||G\psi[a,b] \stackrel{def}{=} \sup_{p \in (a,b)} \left[ \frac{|f|_p}{\psi(p)} \right].
 \end{equation}

 \vspace{4mm}

 \ The function $ \  \psi = \psi(p) = \psi[a,b](p) \  $ is named as  the {\it  generating function } for this space. \par

  \ If for instance

$$
  \psi(p) = \psi^{(r)}(p) = 1, \ p = r;  \  \psi^{(r)}(p) = +\infty,   \ p \ne r,
$$
 where $ \ r = \const \in [1,\infty),  \ C/\infty := 0, \ C \in R, \ $ (an extremal case), then the correspondent
 $ \  G\psi^{(r)}(p)  \  $ space coincides  with the classical Lebesgue - Riesz space $ \ L_r = L_r(Z,\mu). \ $ \par

\vspace{4mm}

 \ Note that the introduced in  \cite{formicagiovamjom2015},  \cite{Greco} etc. norms

$$
||f||_{b, \theta,G} \stackrel{def}{=} \sup_{0 < \epsilon \le b - 1} \left[ \ \epsilon^{\theta/(b-\epsilon)} |f|_{b - \epsilon} \ \right], \ \theta \ge 0
$$
quite coincides with appropriate ones  up to norm equivalence  $ \ ||f||G\psi^{(b; \beta)}, \ \beta = \theta/b. \ $ \par

\vspace{4mm}

 \ Let $ \ f: Z \to R \ $ be certain measurable function such that

$$
\exists a,b, \ 1 \le a < b \le \infty  \ \Rightarrow \forall p \in (a,b) \ |f|_p < \infty.
$$
 \ The so - called {\it natural function} $ \ \psi^{(f)}(p), \ p \in (a,b)  \ $ for this function $ \ f \ $ is defined as follows

$$
\psi^{(f)}(p) \stackrel{def}{=} |f|_p.
$$

 \ Obviously,

$$
||f||G\psi^{(f)} = 1.
$$

\vspace{4mm}

 \ These spaces are investigated in many works, e.g. in
 \cite{Fiorenza1},   \cite{Fiorenza3}, \cite{Fiorenza4},   \cite{Iwaniec1}, \cite{Iwaniec2},
\cite{KozOs1985}, \cite{KozOsSir2017}, \cite{LiflOstSir},   \cite{Ostrovsky1}  - \cite{Ostrovsky5}
etc. They are applied for example in the theory of Partial Differential Equations
\cite{Fiorenza3}, \cite{Fiorenza4}, in the theory of Probability  \cite{Ermakov},\cite{Ostrovsky3}  - \cite{Ostrovsky5}, in Statistics
\cite{Ostrovsky1}, chapter 5, theory of random fields,  \cite{KozOs1985}, \cite{KozOsSir2017}, \cite{Ostrovsky1},  \cite{Ostrovsky4},
in the Functional Analysis \ \cite{Ostrovsky1}, \cite{Ostrovsky2}, \cite{Ostrovsky4} and so one. \par
 \  These spaces are rearrangement invariant (r.i.) Banach functional spaces; its fundamental function  is considered in  \cite{Ostrovsky4}. They do
  not coincides  in general case with the classical rearrangement invariant  spaces: Orlicz, Lorentz, Marcinkiewicz  etc., see \cite{LiflOstSir}, \cite{Ostrovsky2}.

\vspace{5mm}

\section{Main result. The case of small values of parameters.}

\vspace{5mm}

 \hspace{4mm}  We consider in this section the case of Grand Lebesgue Spaces $ \ G\psi_{a,b}, \ $ where $ \ 1 < a < b \le 2. \ $ As before, the measure $ \ \mu \ $
is presumed to be atomless. \par
 \ Introduce the following auxiliary function

$$
\kappa[G\psi[a,b]](u):= \inf_{p \in (a,b)} \left\{ \ \frac{||u||^{2}_p}{\psi^2(p)}   \ \right\}, \ u \in G\psi[a,b].
$$

\ Evidently, this definition is correct also for the arbitrary elements from the  space $ \ L_p: \ u \in L_p, \  1 \le p < \infty, \ $
but we will apply this notion only for the suitable GLS. \par

\vspace{5mm}

 \ {\bf Theorem  2.1.}  In this case, i.e. when $ \ 1 < a < b \le 2, \ $ the  weak characteristic of convexity (WCOC) for
 the space $ \  G\psi = G\psi[a,b] \ $ allows a following lower estimate:

\vspace{4mm}

\begin{equation} \label{small ab}
\Delta[G\psi[a,b]](x - y) \ge \frac{a-1}{4} \cdot   \kappa[G\psi[a,b]](x - y), \ x,y \in B[G\psi[a,b]],
\end{equation}

so that

$$
||x + y||G\psi  \le 2 - \frac{a-1}{4} \kappa[G\psi[a,b]](x - y),  \ x,y \in B[G\psi[a,b]].
$$

\vspace{4mm}

 \ {\bf Proof.}   Let $ \  x,y \in B[G\psi[a,b]]. \  $  We have from the direct definition of the
 norm in GLS

$$
\frac{||x||_p}{\psi(p)} \le 1, \hspace{3mm} \frac{||y||_p}{\psi(p)} \le 1, \ p \in (a,b).
$$

 \ It follows on the basis of definition 1.1 being applied to the space $ \ L_p \ $

$$
\left| \ \left| \ \frac{x}{\psi(p)} + \frac{y}{\psi(p)} \  \right| \ \right|_p \le 2 - 2  \delta_{L_p} \left( \ \frac{||x - y||_p}{\psi(p)}  \ \right).
$$

 \ We derive taking supremum over $ \ p \ $ using the estimate

\begin{equation} \label{General p}
||x + y||G\psi[a,b] \le 2 - 2 \inf_{p \in (a,b)} \delta_{L_p} \left( \  \frac{||x - y||_p}{\psi(p)}  \ \right).
\end{equation}

 \vspace{3mm}

 \ Further, we will use in particular the relation  (\ref{small p}); from one and (\ref{General p})  it follows

\begin{equation} \label{small p  a1}
\delta_p(\epsilon) \ge \frac{a-1}{8} \ \epsilon^2, \ \epsilon \in [0,2].
\end{equation}

 \ Therefore

$$
\frac{||x + y||_p}{\psi(p)} \le 2 - \frac{a-1}{4} \cdot \frac{||x - y||_p^2}{\psi^2(p)}.
$$

 \ It remains to take the supremum over $ \ p \in (a,b), \ $ using the direct definition of the norm in the Grand Lebesgue Spaces:

$$
||x + y||G\psi  = \sup_{p \in (a,b)} \left\{ \ \frac{||x +y||_p}{\psi(p)} \ \right\} \le
$$

$$
 \le 2 -   \frac{a-1}{4} \ \inf_{p \in (a,b)} \left\{ \ \frac{||x-y||_p^2}{\psi^2(p)} \ \right\} = 2 - \frac{a-1}{4} \kappa[G\psi[a,b]](x - y),
$$
Q.E.D. \par

\vspace{5mm}

\section{Main result. The case of great values of parameters.}

\vspace{5mm}

 \hspace{4mm} We consider in this section the opposite  case of Grand Lebesgue Spaces $ \ G\psi_{a,b}, \ $ where $ \ 2 < a < b  < \infty. \ $
The measure $ \ \mu \ $ is again presumed to be atomless. \par

 \ Introduce the following auxiliary function

$$
\theta[G\psi[a,b]](u):= \inf_{p \in (a,b)} \left\{ \ \frac{||u||^{p}_p}{ p \cdot 2^p \cdot \psi^p(p)}   \ \right\}, \ u \in G\psi[a,b].
$$
 \ Obviously, the last definition is correct also for the space $ \ L_p: \ u \in L_p, \  1 \le p < \infty. \ $ \par

\vspace{5mm}

 \ {\bf Theorem  3.1.}  In this case the  weak characteristic of convexity (WCOC) for the space $ \  G\psi = G\psi[a,b], \ $ where $ \ 2 < a < b < \infty \ $ obeys
a following lower estimate:

\vspace{4mm}

\begin{equation} \label{small ab}
\Delta[G\psi[a,b]](x - y) \ge \theta[G\psi[a,b]](x - y),
\end{equation}

so that

$$
||x + y||G\psi  \le 2 -  \theta[G\psi[a,b]](x - y),  \ x,y \in B[G\psi[a,b]].
$$

\vspace{4mm}

 \ {\bf Proof} \  is quite alike as before in the previous section, as well. We will use the relation
  (\ref{great p appr}) for the values  $ \  x,y \in B[G\psi[a,b]]: \ $

$$
\frac{||x + y||_p}{\psi(p)} \le 2 - \frac{||x-y||_p^p}{p  \cdot 2^p \cdot \psi^p(p)}, \ x,y \in B[G\psi[a,b]] .
$$

 \ It remains to take the supremum over $ \ p \in (a,b), \ $ using the direct definition of the norm in the Grand Lebesgue Spaces:

$$
||x + y||G\psi  = \sup_{p \in (a,b)} \left\{ \ \frac{||x +y||_p}{\psi(p)} \ \right\} \le
$$

$$
 \le 2 - \inf_{p \in (a,b)} \left\{ \ \frac{||x-y||_p^{p}}{ p \cdot 2^p \cdot \psi^p(p)} \ \right\} = 2 - \theta[G\psi[a,b]](x - y),
$$
Q.E.D. \par

\vspace{5mm}

\section{ Some examples.}

\vspace{4mm}

\hspace{3mm} {\bf Example 1.} Suppose in addition to the of theorem 2.1., i.e.
 to the case of Grand Lebesgue Spaces $ \ G\psi_{a,b}, \ $ where $ \ 1 < a < b \le 2, \ $ that the measure $ \ \mu \ $
is bounded: $ \ \mu(Z) = 1. \ $  Then one can apply the famous  Liyapunov's   inequality: $ \ ||x||_p \ge ||x||_a, \ p \in [a,b].  \ $ \par

 \ Assume yet that the generating function $ \ \psi(\cdot) \ $ is upper bounded: $ \ \psi(p) \le d = \const \in (0,\infty). \ $  It follows from the
proposition of theorem 2.1

$$
||x + y||G\psi[a,b] \le 2 - \frac{a-1}{4} \ \frac{||x-y||_a^2}{d^2}, \ x,y \in B[G\psi[a,b]].
$$

\vspace{4mm}

 {\bf Example 2.} We retain the restrictions $ \ \mu(Z) = 1 \ $ and  $ \ \psi(p) \le d = \const \in (0,\infty); \ $
 but let here $ \ 2 \le a < b < \infty. \ $ Then
$ \ ||x||_p \le ||x||_b, \ p \le b;  \ $ and we deduce under conditions and assertions of theorem 2.1

$$
||x + y||G\psi[a,b] \le 2 -  \frac{||x-y||_a^a}{b \cdot 2^b \ d^b}, \ x,y \in B[G\psi[a,b]].
$$

\vspace{5mm}

\section{Concluding remarks.}

\vspace{4mm}

 \hspace{3mm}  It is interest in our opinion to study the Modulus Of Convexity  for other Grand Lebesgue Spaces, in
 particular, which are mentioned in the first section; especially in the case when $ \ b = \infty. \ $ \par
 \ Open question: are the GLS $ \ G\psi[a,\infty], \ $ for instance, the Subgaussian Spaces with

 $$
 \psi(p) = \sqrt{p}, \ p \in (1,\infty)
 $$
modulative convex? \par

\vspace{5mm}

\vspace{0.5cm} \emph{Acknowledgement.} {\footnotesize The first
author has been partially supported by the Gruppo Nazionale per
l'Analisi Matematica, la Probabilit\`a e le loro Applicazioni
(GNAMPA) of the Istituto Nazionale di Alta Matematica (INdAM) and by
Universit\`a degli Studi di Napoli Parthenope through the project
\lq\lq sostegno alla Ricerca individuale\rq\rq (triennio 2015 - 2017)}.\par

\end{document}